\documentclass[12pt]{article}
\usepackage[cp1251]{inputenc}
\usepackage[russian]{babel}
\usepackage{amssymb}
\usepackage{amsmath}
\begin{document}

\begin{center}
\textbf{\large On uniform convergence of Fourier series}
\end{center}

\begin{center}
Vladimir Lebedev
\end{center}

\begin{quotation}
{\small \textsc{Abstract.} We consider the space $U(\mathbb T)$
of all continuous functions on the circle $\mathbb T$ with
uniformly convergent Fourier series. We show that if $\varphi:
\mathbb T\rightarrow\mathbb T$ is a continuous piecewise linear
but not linear map, then $\|e^{in\varphi}\|_{U(\mathbb
T)}\simeq\log n$.

  References: 5 items.


  }

\end{quotation}

\quad

   Given any integrable function
$f$ on the circle $\mathbb T=\mathbb R /2\pi \mathbb Z$ (where
$\mathbb R$ is the real line and $\mathbb Z$ is the set of
integers), consider Fourier series
$$
f(t)\sim \sum_{k\in \mathbb Z}\widehat f(k) e^{ikt}, \qquad
\widehat{f}(k)=\frac{1}{2\pi}\int_{\mathbb T} f(t)e^{-ikt} dt.
$$

   Let $C(\mathbb T)$ be the space of continuous functions
$f$ on $\mathbb T$ with the usual norm $\|f\|_{C(\mathbb
T)}=\sup_{t\in\mathbb T}|f(t)|$.

   Let $U(\mathbb T)$ be the space of all functions
$f\in C(\mathbb T)$ whose Fourier series converges uniformly, that
is $\|S_N(f)-f\|_{C(\mathbb T)}\rightarrow 0$ as
$N\rightarrow\infty$, where $S_N(f)$ stands for the $N$ -th
partial sum of the Fourier series of $f$:
$$
S_N(f)(t)=\sum_{|k|\leq N}\widehat f(k) e^{ikt}.
$$
Endowed with the natural norm
$$
\|f\|_{U(\mathbb T)}=\sup_N\|S_N(f)\|_{C(\mathbb T)},
$$
the space $U(\mathbb T)$ is a Banach space.

  Consider also the space $A(\mathbb T)$ of functions
$f\in C(\mathbb T)$ with absolutely convergent Fourier series. We
put
$$
\|f\|_{A(\mathbb T)}=\|\widehat{f}\|_{l^1}=
\sum_{k\in\mathbb Z} |\widehat{f}(k)|.
$$
The space $A(\mathbb T)$ is a Banach space. We have $A(\mathbb
T)\subset U(\mathbb T)$ and obviously $\|\cdot\|_{U(\mathbb
T)}\le\|\cdot\|_{A(\mathbb T)}$.

  Let $\varphi$ be a continuous map of the circle
$\mathbb T$ into itself, i.e. a continuous function $\varphi:
\mathbb R\rightarrow\mathbb R$ satisfying
$\varphi(t+2\pi)=\varphi(t)~(\mathrm{mod}\,2\pi), ~t\in\mathbb
R$. According to the known Beurling--Helson theorem [1] (see
also [2], [3]), if $\|e^{in\varphi}\|_{A(\mathbb T)}=O(1),
~n\rightarrow\infty$, then $\varphi$ is linear. At the same
time it is known that there exist nontrivial maps $\varphi :
\mathbb T\rightarrow\mathbb T$ such that
$\|e^{in\varphi}\|_{U(\mathbb T)}=O(1)$. For a survey of
certain results on the maps of the circle and the spaces
$A(\mathbb T), ~U(\mathbb T)$ see [3], [4]. For more recent
results see [5].

  The following assertion is due to J.-P. Kahane
[2, Ch. IV]. If a nonlinear continuous map $\varphi : \mathbb
T\rightarrow\mathbb T$ is piecewise linear (which means that $[0,
2\pi]$ is a finite union of intervals such that $\varphi$ is
linear on each of them), then $\|e^{in\varphi}\|_{A(\mathbb
T)}\simeq\log |n|$. (The sign $\simeq$ means that for all
sufficiently large $n\in\mathbb Z$ the ratio of the corresponding
quantities is contained between two positive constants.) Here we
shall obtain a similar assertion for the space $U(\mathbb T)$.

\quad

\textbf{Theorem.} \emph {Let $\varphi$ be a piecewise linear but
not linear continuous map of the circle $\mathbb T$ into itself.
Then} $ \|e^{in\varphi}\|_{U(\mathbb T)}\simeq\log |n|,
~n\in\mathbb Z.$

\quad

  In particular, this theorem implies that, generally speaking,
nontrivial piecewise linear changes of variable destroy the
uniform convergence of Fourier series. Moreover they do not act
from $A(\mathbb T)$ to $U(\mathbb T)$. Indeed, assuming that
for each function $f\in A(\mathbb T)$ the superposition
$f\circ\varphi$ belongs to $U(\mathbb T)$, we would have
$\|e^{in\varphi}\|_{U(\mathbb T)}=O(1)$ (it suffices to apply
the closed graph theorem to the operator $f\rightarrow
f\circ\varphi$).

  To prove the theorem we only have to prove
$\log |n|$ lower bound, the upper bound
$\|e^{in\varphi}\|_{U(\mathbb T)}=O(\log |n|)$ follows from
inequality $\|\cdot\|_{U(\mathbb T)}\leq \|\cdot\|_{A(\mathbb
T)}$ and the above result of Kahane.

  We shall need the following simple lemma that perhaps is of interest
in itself.

\quad

\textbf{Lemma.} \emph {Let $m\in C(\mathbb T)$ be a function such
that
$$
\|m\|_*=\sum_{n\in\mathbb Z}|\widehat{m}(n)|\log(|n|+2)<\infty.
$$
Then for each function $f\in U(\mathbb T)$ we have $mf\in
U(\mathbb T)$ and
$$
\|mf\|_{U(\mathbb T)}\leq c\|m\|_*\,\|f\|_{U(\mathbb T)},
$$
where the constant $c>0$ is independent of $f$ and $m$.}

\quad

\emph{Proof.} Let $f\in U(\mathbb T)$. For $n\in \mathbb Z$ we
put $e_n(t)=e^{int}$. Let $n>0$. Then
$$
S_N(e_n f)=e_nS_{N+n}(f)+e_N\widehat{f}(N-n)-e_{N+n}S_n(e_{-N}f).
\eqno(1)
$$
It is easy to see that this relation implies the inclusion
$e_nf\in U(\mathbb T)$.

  For each function $g\in U(\mathbb T)$
we have $\|g\|_{C(\mathbb T)}\leq \|g\|_{U(\mathbb T)}$. At the
same time (as is well known) for each function $g\in C(\mathbb T)$
we have
$$
\|S_n(g)\|_{C(\mathbb
T)}\leq c\|g\|_{C(\mathbb T)}\log (n+2)
$$
with a constant $c>0$ independent of $n$ and $g$. So (1) yields
$$
\|S_N(e_nf)\|_{C(\mathbb T)}\leq \|S_{N+n}(f)\|_{C(\mathbb
T)}+\|f\|_{C(\mathbb T)}+\|S_n(e_{-N}f)\|_{C(\mathbb T)}
$$
$$
\leq\|f\|_{U(\mathbb T)}+\|f\|_{C(\mathbb T)}+
c\|f\|_{C(\mathbb T)}\log (n+2)\leq c_1\|f\|_{U(\mathbb T)}\log (n+2).
$$
Thus, $\|e_nf\|_{U(\mathbb T)}\leq c_1\|f\|_{U(\mathbb T)}\log
(|n|+2)$. The same relation holds for $n<0$ (complex
conjugation does not affect the norm of a function in
$U(\mathbb T)$). The assertion of the lemma immediately
follows.

\quad

   \emph{Proof of the theorem.} In the standard way we identify
functions on $\mathbb T$ with functions on the interval $[-\pi,
\pi]$. For $v\in\mathbb R$ define the functions $e_v$ on $[-\pi,
\pi]$ by $e_v(t)=e^{ivt}$. For an arbitrary interval $I\subseteq
[-\pi, \pi]$ let $1_I$ denotes its characteristic function:
$1_I(t)=1$ for $t\in I$, $1_I(t)=0$ for $t\in [-\pi, \pi]\setminus
I$. For $0<\varepsilon<\pi$ let $\Delta_\varepsilon$ be the
``triangle'' function supported on the interval $(-\varepsilon,
\varepsilon)$, that is a function on $[-\pi, \pi]$ defined by
$\Delta_\varepsilon (t)=\max (0, 1-|t|/\varepsilon)$.

   Let $t_0$ be a point such that $\varphi$ is linear in
some left half-neighborhood of $t_0$ and is linear in some
right half-neighborhood of $t_0$, but is not linear in any its
neighborhood. Replacing (if necessarily) the function
$\varphi(t)$ by $\varphi(t+t_0)-\varphi(t_0)$ we can assume
that $t_0=0$ and $\varphi(t_0)=0$; thus we can assume that for
a certain $\varepsilon, ~0<\varepsilon<\pi,$ we have
$\varphi(t)=\alpha t$ for $t\in(-\varepsilon, 0]$ and
$\varphi(t)=\beta t$ for $t\in[0, \varepsilon)$, where
$\alpha\neq \beta$.

  Direct calculation yields for $k\neq n\alpha, ~n\beta$
$$
\widehat{\Delta_\varepsilon e^{in\varphi}}(k)=\frac{1}{2\pi
i}\bigg(\frac{1}{n\alpha -k}-\frac{1}{n\beta -k}\bigg)-
$$
$$
-\frac{1}{i\varepsilon}\bigg(\frac{1}{n\alpha-k}\widehat{e_{n\alpha}
1_{(-\varepsilon, 0)}}(k)- \frac{1}{n\beta-k}\widehat{e_{n\beta}
1_{(0, \varepsilon)}}(k)\bigg).
\eqno(2)
$$

   For
$\lambda\in\mathbb R$ let
$$
Q(\lambda)=\sum_{k\in\mathbb Z : ~|k-\lambda|\geq
1}\frac{1}{(k-\lambda)^2}.
$$
It is easy to verify that
$$
Q(\lambda)\leq 4, \qquad \lambda\in\mathbb R.
\eqno (3)
$$

  We shall show first that for $n\in\mathbb Z, ~n\neq 0,$ we have
$$
\|\Delta_\varepsilon e^{in\varphi}\|_{U(\mathbb T)}\geq \frac{1}{2\pi}\log
|n|+c(\varphi),
\eqno (4)
$$
where $c(\varphi)$ is independent of $n$.

  If $g\in U(\mathbb T)$ then the function $g(-t)$
and the function $\overline{g(t)}$ (obtained by complex
conjugation) belong to $U(\mathbb T)$ and have the same norm as
that of $g$. So in the proof of estimate (4) we can assume that
$\alpha>0$ and consider only the following three cases: 1)
$|\beta|>\alpha$; 2) $\beta=-\alpha$; 3) $\beta=0$.

   We can also assume that $n$ is positive and is so large that
$n\alpha\geq 2$.

   In what follows we put $N=[n\alpha]-1$, where $[x]$
stands for the integer part of a number $x$.

  \emph{Case} 1). We have
$$
\bigg|\sum_{|k|\leq N}\frac{1}{n\alpha-k}\bigg|=\sum_{|k|\leq
N}\frac{1}{n\alpha-k}\geq \sum_{|k|\leq
N}\frac{1}{N+2-k}
$$
$$
=\frac{1}{2}+\frac{1}{3}+\ldots +\frac{1}{2N+2}\geq
\log(N+1)\geq\log\frac{n\alpha}{2}.
\eqno (5)
$$
At the same time for $|k|\leq N$ we have $|n\beta-k|\geq
|n\beta|-|n\alpha|=n(|\beta|-\alpha)$, so
$$
\bigg|\sum_{|k|\leq N}\frac{1}{n\beta-k}\bigg|\leq
\frac{2N+1}{n(|\beta|-\alpha)}\leq\frac{3n\alpha}{n(|\beta|-\alpha)}=
\frac{3\alpha}{|\beta|-\alpha}. \eqno (6)
$$

  Note now that using the Cauchy inequality and the Parseval
identity we obtain (see (3))
$$
\bigg|\sum_{|k|\leq N}\frac{1}{n\alpha-k}\widehat{e_{n\alpha}
1_{(-\varepsilon, 0)}}(k)\bigg|\leq (Q(n\alpha))^{1/2}\|e_{n\alpha}
1_{(-\varepsilon, 0)}\|_{L^2(\mathbb T)}\leq 2\varepsilon^{1/2}, \eqno
(7)
$$
and similarly
$$
\bigg|\sum_{|k|\leq N}\frac{1}{n\beta-k}\widehat{e_{n\beta} 1_{(0,
\varepsilon)}}(k)\bigg|\leq (Q(n\beta))^{1/2}\|e_{n\beta} 1_{(0,
\varepsilon)}\|_{L^2(\mathbb T)}\leq 2\varepsilon^{1/2}. \eqno (8)
$$

   Together relations (5)--(8) imply (see (2))
$$
|S_N(\Delta_\varepsilon e^{in\varphi})(0)|=\bigg|\sum_{|k|\leq N}
\widehat{\Delta_\varepsilon e^{in\varphi}}(k)\bigg|\geq
\frac{1}{2\pi}\bigg(\log\frac{n\alpha}{2}-\frac{3\alpha}
{|\beta|-\alpha}\bigg)-4\varepsilon^{-1/2},
$$
and we obtain (4).

   \emph{Case} 2). We have
$$
\bigg|\sum_{|k|\leq
N}\bigg(\frac{1}{n\alpha-k}-\frac{1}{n\beta-k}\bigg)\bigg|=
\bigg|\sum_{|k|\leq
N}\bigg(\frac{1}{n\alpha-k}+\frac{1}{n\alpha+k}\bigg)\bigg|
$$
$$
=\bigg|2\sum_{|k|\leq N}\frac{1}{n\alpha-k}\bigg|\geq
\log\frac{n\alpha}{2}.
$$
Together with estimates (7), (8), which are true in Case 2), this
estimate yields
$$
|S_N(\Delta_\varepsilon e^{in\varphi})(0)|\geq
\frac{1}{2\pi}\log\frac{n\alpha}{2}-4\varepsilon^{-1/2},
$$
and we obtain (4) again.

  \emph{Case} 3). We have
$$
\bigg|\sum_{1\leq |k|\leq
N}\bigg(\frac{1}{n\alpha-k}-\frac{1}{n\beta-k}\bigg)\bigg|=
\bigg|\sum_{1\leq |k|\leq
N}\bigg(\frac{1}{n\alpha-k}-\frac{1}{-k}\bigg)\bigg|
$$
$$
=\bigg|\sum_{1\leq |k|\leq N}\frac{1}{n\alpha-k}\bigg|\geq
\log\frac{n\alpha}{2}-1.
$$
Note that estimates (7), (8) are true in Case 3) if we replace the
range $|k|\leq N$ in the sums by $1\leq |k|\leq N$. Thus we see
that
$$
\bigg|\sum_{1\leq |k|\leq N} \widehat{\Delta_\varepsilon
e^{in\varphi}}(k)\bigg|\geq
\frac{1}{2\pi}\bigg(\log\frac{n\alpha}{2}-1\bigg)-4\varepsilon^{-1/2},
$$
and since $|\widehat{\Delta_\varepsilon e^{in\varphi}}(0)|\leq 1$,
we obtain
$$
|S_N(\Delta_\varepsilon e^{in\varphi})(0)|\geq
\frac{1}{2\pi}\bigg(\log\frac{n\alpha}{2}-1\bigg)-4\varepsilon^{-1/2}-1.
$$
Estimate (4) is proved.

   Note now that
$\widehat{\Delta_\varepsilon}(k)=O(1/|k|^2)$ as
$|k|\rightarrow\infty$, so
$$
\|\Delta_\varepsilon\|_*=\sum_{k\in\mathbb
Z}|\widehat{\Delta_\varepsilon}(k)|\log(|k|+2)=M(\varepsilon)<\infty,
$$
and from (4), using Lemma, we obtain
$$
c(\varphi)+\frac{1}{2\pi}\log |n|\leq \|\Delta_\varepsilon
e^{in\varphi}\|_{U(\mathbb T)}\leq
c M(\varepsilon) \|e^{in\varphi}\|_{U(\mathbb T)}.
$$
The theorem is proved.

\quad

\emph{Remark.} It would be interesting to describe pointwise
multipliers of the space $U(\mathbb T)$, i.e. the continuous
functions $m$ on $\mathbb T$ such that $mf\in U(\mathbb T)$
whenever $f\in U(\mathbb T)$. According to the lemma obtained
above the condition
$$
\sum_{k\in\mathbb Z}|\widehat{m}(k)|\log(|k|+2)<\infty
$$
is a sufficient condition for a function $m$ to be a
multiplier. The author does not know if this condition is
necessary. The indicated condition can not be replaced by the
weaker condition $m\in A(\mathbb T)$ (see [2, Ch. I, \S ~6]).

\begin{center}
\textsc{References}
\end{center}

\begin{enumerate}

\item A. Beurling, H. Helson, ``Fourier-Stieltjes transforms
    with bounded powers'', \emph{Math. Scand.,} \textbf{1}
    (1953), 120-126.

\item J.-P. Kahane, \emph{S\'eries de Fourier absolument
    convergantes}, Springer-Verlag, Berlin--Heidelberg--New
    York, 1970.

\item J.-P. Kahane, ``Quatre le\c cons sur les
    hom\'eomorphismes du circle et les s\'eries de Fourier'',
    in: \emph{Topics in Modern Harmonic Analysis,} Vol. II,
    Ist. Naz. Alta Mat. Francesco Severi, Roma, 1983, 955-990.

\item A. M. Olevskii, ``Modifications of functions and Fourier
    series'', \emph{Russian Math. Surveys}, \textbf{40}:3
    (1985), 181-224.

\item  V. V. Lebedev, ``Quantitative estimates in Beurling
    --Helson type theorems'', \emph{Sbornik: Mathematics},
    \textbf{201}:12 (2010), 1811-1836.

\end{enumerate}

\quad

\noindent Dept. of Mathematical Analysis\\
Moscow State Institute of Electronics\\
and Mathematics (Technical University)\\
E-mail address: \emph {lebedevhome@gmail.com}

\end{document}